\documentclass{amsart}
\usepackage{amsmath, amsfonts, amssymb}
\usepackage{xypic}
\usepackage{graphicx}

\def \Xh {{{\widehat{X}}}}
\newtheorem{thm}{Theorem}[section] 
\newtheorem{lemma}[thm]{Lemma}
\newtheorem{prop}[thm]{Proposition}
\newtheorem{cor}[thm]{Corollary}

\theoremstyle{definition}
\newtheorem{remark}[thm]{Remark}
\newtheorem{example}[thm]{Example}
\newtheorem{defn}[thm]{Definition}

\DeclareMathOperator{\Hom}{Hom}

\DeclareMathOperator{\codim}{codim}

\newcommand{\bP}{{\mathbb P}}

\newcommand{\G}{{\mathbb G}}

\renewcommand{\cL}{{\mathcal L}}

\newcommand{\cO}{{\mathcal O}}

\newcommand{\ignore}[1]{}

\def \im {{{\rm Im}}}

\def \diag {{{\rm diag}}}

\newcommand{\scal}[1]{\langle #1 \rangle}

\def \cP {{{\mathcal{P}}}}
\def \cH {{{\mathcal{H}}}}
\def \cW {{{\mathcal{W}}}}
\def \cT {{{\mathcal{T}}}}
\def \cG {{{\mathcal{G}}}}
\def \cX {{{\mathcal{X}}}}
\def \cY {{{\mathcal{Y}}}}
\def \cU {{{\mathcal{U}}}}
\def \cB {{{\mathcal{B}}}}

\def \Xt {{{\widetilde{X}}}}
\def \It {{{\widetilde{I}}}}
\def \Jt {{{\widetilde{J}}}}
\def \Kt {{{\widetilde{K}}}}

\def \xt {{{\widetilde{x}}}}
\def \yt {{{\widetilde{y}}}}
\def \Ht {{{\widetilde{H}}}}
\def \Kt {{{\widetilde{K}}}}
\def \Rt {{{\widetilde{R}}}}

\def \Ot {{{\widetilde{\Omega}}}}
\def \Dt {{{\widetilde{D}}}}
\def \Xh {{{\widehat{X}}}}
\def \Gh {{{\widehat{G}}}}
\def \Hh {{{\widehat{H}}}}

\def \Xo {{\mathring{X}}}
\def \Pio {{\mathring{\Pi}}}
\def \Ro {{\mathring{R}}}
\def \PRo {{\mathring{PR}}}
\def \Rto {{\mathring{\Rt}}}
\def \Xto {{\mathring{\Xt}}}
\def \Yo {{\mathring{Y}}}

\def \Pio {{\mathring{\Pi}}}
\def \a {{{\alpha}}}

\begin{document}

\title{Split subvarieties of group embeddings}


\author{Nicolas Perrin}
\address{Mathematisches Institut, Heinrich-Heine-Universit{\"a}t,
D-40204 D{\"u}sseldorf, Germany}
\email{perrin@math.uni-duesseldorf.de}



\begin{abstract}
Let $G$ be a connected reductive group and $X$ an equivariant compactifiction of
$G$. In $X$, we study generalised and opposite generalised Schubert
varieties, their intersections called generalised Richardson varieties
and projected generalised Richardson varieties. Any complete
$G$-embedding has a canonical Frobenius splitting and we prove that
the compatibly split subvarieties are the generalised projected
Richardson varieties extending a result of Knutson, Lam and Speyer to
the situation.
\end{abstract}

\maketitle

\markboth{N.~PERRIN}{Split subvarieties of group embeddings}


\section{Introduction}\label{sec:intro}

Let $G$ be a connected reductive group over a field $k$ of positive characteristic $p$. A $G$-embedding is a normal $G\times G$-variety $X$ together with a $G\times G$-equivariant open embedding of $G$ in $X$.

$G$-embeddings share many of the nice properties of rational
projective homogeneous spaces. For example, any $G$-embedding has
a cellular decomposition defined by $B \times B$ and $B^- \times
B^-$-orbits for $B$ and $B^-$ opposite Borel subgroups of $G$. We call
these cells and their closures generalised (and opposite generalised)
Schubert cells and varieties. They have many properties of the
classical Schubert cells and varieties: normality, Cohen-Macaulay
property (see for example \cite{BP,BT,he-thomsen} for more
details). We study these varieties and their intersections, that we
call generalised Richardson varieties, as well as the images of these
varieties under morphisms of $G$-embeddings. 

The existence of Frobenius splittings is another instance of the
common features between projective rational homogeneous spaces and
$G$-embeddings. Frobenius splittings were first introduced by Mehta
and Ramanathan in \cite{MR} for projective rational homogeneous spaces
to prove cohomology vanishing results and regularity properties of
Schubert varieties. Using this technique, Rittatore \cite{rittatore}
obtained regularity results for all $G$-embeddings, in particular the
Cohen-Macaulay property. Brion and Polo \cite{BP}, Brion and Thomsen
\cite{BT} and He and Thomsen \cite{he-thomsen} also obtained
regularity results for $B\times B$-orbit closures in group embeddings.

For rational projective homogeneous spaces, Knutson, Lam and Speyer
\cite{KLS} proved that in $G/P$ (with $P$ a parabolic subgroup
containing $B$) the projections of Richardson varieties are all the
compatibly split subvarieties for the unique $B$-canonical
splitting. For $X$ a complete $G$-embedding, $X$ has a unique
Frobenius splitting $\phi$ compatibly splitting the $G\times G$,
$B\times B$ and $B^-\times B^-$ divisors (see Proposition
\ref{thm-HT}). We introduce projected generalised Richardson varieties
(see Definition \ref{defi-projected}) and prove the following result. 

\begin{thm}
The projected generalised Richardson varieties are the
$\phi$-compatibly split subvarieties. 
\end{thm}

The fact that projected generalised Richardson varieties are
compatibly split follows from results of He and Thomsen
\cite{he-thomsen}. We use techniques of Knutson, Lam and Speyer
\cite{KLS} to prove that these varieties are the only compatibly split
subvarieties. On the way we prove several results on generalised
Schubert varieties, generalised Richardson varietes and projected
generalised Richardson varieties: the later are normal and form
a stratification of $X$. It is also interesting to note that
not all the properties of Schubert varieties extend to
$G$-embeddings. In particular, in general the intersection of two
opposite generalised Schubert varieties is neither irreducible nor
equidimensional (see Example \ref{exam-1}). 

\subsection*{Acknowledgements} 

I would like to thank Michel Brion for many suggestions and comments. Especially the idea of using of algebraic monoids in the proof of Proposition \ref{prop-stab} is due to him. I also thank the participants of the seminar on wonderful compactifications in Essen which was the starting point of this project.

\subsection*{Notation}

We work over an algebraically closed field $k$ of positive characteristic $p$. Varieties are reduced, separated, connected schemes of finite type over $k$. 

Let $G$ be a reductive group over $k$ and let $T$ be a maximal torus
of $G$. Denote by $W=N_G(T)/T$ the Weyl group of $T$ and by $\Phi$ the
root system associated to $(G,T)$. Let $B$ be a Borel subgroup of $G$
containing $T$. Denote by $\Delta$ the set of simple roots induced by
$B$ and by $\Phi^+$ the set of positive roots. For $J$ a subset of
$\Delta$, denote by $P_J$ the parabolic subgroup containing $B$ with
$\Delta_P = J$ where for $P$ a parabolic subgroup containing $B$,
$\Delta_P$ is the set of simple roots of the Levi factor of $P$
containing $T$. Denote by $P_J^-$ the opposite parabolic subgroup and
by $L_J$ the Levi subgroup containing $T$ of both $P_J$ and
$P_J^-$. Write $Z_J$ for the center of $L_J$. We write $U_J$ and
$U_J^-$ for the unipotent radicals of $P_J$ and $P_J^-$. We write
$W_J$ for the Weyl group of $P_J$ and $W^J$ for  the set of minimal
length representatives of $W/W_J$. Recall that there exists for $u\in
W$ a unique length additive decomposition $u=u^Ju_J$ with $u^J\in W^J$
and $u_J\in W_J$. Denote by $B_J$ the intersection $B\cap L_J$ and by
$B_J^-$ the intersection $B^-\cap L_J$. We write $w_0$ for the longest
element in $W$. For $L$ a group, we denote by $(L,L)$ its derived
group and $\diag(L)$ the diagonal embedding of $L$ in $L \times L$.

\section{$G$-embeddings}

\subsection{Toroidal $G$-embeddings}

Consider the $G\times G$-action on $G$ given by $(g_1,g_2)\cdot
g=g_1gg_2^{-1}$. A $G$-embedding is a normal $G\times G$-variety $X$
together with an open equivariant embedding $G\to X$. A morphism of
$G$-embeddings is a $G \times G$-equivariant morphism between
$G$-embeddings extending the identity on $G$. These varieties
are special cases of spherical varieties. We refer to
\cite{knop,survey} for an overview on the geometry of spherical
varieties. 

\begin{defn}
\label{def-toro}
A $G$-embedding $X$ is called toroidal if any $B\times B$-stable
divisor in $X$ containing a $G\times G$-orbit is $G\times G$-stable. A
$G$-embedding $X$ is called simple if $X$ has a unique closed $G\times
G$-orbit in $X$.
\end{defn}

\subsection{Description of $G \times G$-orbits}

Let $G$ be a reductive group, $X$ be a $G$-embedding and $(G \times G) \cdot x$
be a $G \times G$-orbit in $X$. We describe the stabiliser $H$ of
$x$. The following result, whose proof is essentially due to Brion,
generalises in positive characteristic a result of Alexeev and Brion
\cite[Proposition 3.1]{BA}

\begin{prop}
\label{prop-stab}
1. There exists a element $x' \in (G \times G) \cdot x$ unique up to $T \times T$-action such that $\diag(T)$ fixes $x'$.

2. Assume that $x = x'$. Then there exists a subset $I$ of $\Delta$, union of two orthogonal subsets $J$ and $K$ such that the subgroup $H$ is conjugate in $T \times T$ to
$(U_I H_J \times U^-_I H_J) \diag L_K$ where $(L_J,L_J) \subset H_J
\subset (L_J,L_J)Z_K$ and $Z_K$ is the center of $L_K$. 

Furthermore if $X$ is toroidal, then $J = \emptyset$.
\end{prop}

\begin{proof}
1. This follows from \cite[Proposition 6.2.3]{BK} for toroidal embeddings. The general result follows from the toroidal case. Note that since $x'$ is unique up to $T \times T$-action, it follows that its stabiliser will be unique up to conjugation in $T \times T$.

2. Using results of Sumihiro \cite{sumihiro} (see also \cite[Theorem
  2.3.1]{survey}), we may assume that $X$ is equivariantly embedded in $\mathbb{P}(V)$ with $V$ a $G$-module. Consider $\Xh$ the affine
cone over $X$ and $\Gh = G \times \G_m$. The stabiliser of the cone
over the orbit $G \cdot x$ is $\Hh \simeq H$. We can thus
assume that $X$ is affine. According to a result of Rittatore (see
\cite[Proposition 1]{rittatore98} the affine $G$-embedding $X$ is an
algebraic monoid. The result is a consequence of the theory of
algebraic monoids. For this theory, we refer to \cite{renner} although
many of the results we use were first proved by Putcha
\cite{putcha1,putcha2}. 

By \cite[Theorem 4.5.(c)]{renner}, any $G \times G$-orbit is the orbit of an
idempotent $e$ contained in the closure of the maximal torus. We may
therefore replace $x$ by $e$. The stabiliser is the subgroup
$H=\{(x,y)\in G\times G\ |\ xey^{-1}=e\}$ of $G\times G$. Let $(x,y)
\in H$, then $xe=ey$ and $exe=ey=xe$. Therefore $x$ lies in
$P(e)=\{x\in G\ |\ xe=exe\}$. In the same way, $eye=xe=ey$ and $y$
lies in $P^-(e)=\{y\in G\ |\ ey=eye\}$. According to \cite[Theorem
  4.5.(a)]{renner}, these groups are opposite parabolic subgroups of
$G$ and their unipotent radicals $R_uP(e)$ and $R_uP^-(e)$ satisy
$R_uP(e)e=\{e\}=eR_uP^-(e)$. In particular we have the inclusions 
$$R_uP(e) \times R_uP^-(e) \subset H \subset P(e) \times P^-(e).$$
Note that the Levi subgroup of both $P(e)$ and $P^-(e)$ is $L(e) =
P(e)\cap P^-(e) = \{ x \in G \ | \ xe=ex \} = C_G(e)$.  

By \cite[Theorem 4.8.(a)]{renner}, the subset $eXe = \{x \in X\ |\ x =
exe\}$ is an algebraic monoid with unit $e$ and unit group
$C_G(e)e=eC_G(e)$. Consider the morphism $p_e : P(e) \times P^-(e) \to
C_G(e)e \times e C_G(e)$ defined by $p_e(x,y)=(xe,ey)$. It is a group
homomorphism: $p_e(xx',yy')=(xx'e,yy'e)=(xex'e,yey'e)$, whose kernel
contains $R_uP(e)\times R_uP^-(e)$. Thus $p_e$ factors through its
restriction to $L(e) \times L(e)$. Since $L(e)$ is reductive, the
morphism $L(e) \to C_G(e)e,\ x\mapsto xe$ is the quotient of a finite cover of $L(e)$ by some semi-simple factors and a subgroup of the centre. 

For $(x,y) \in H$, we have $p_e(x,y)=(xe,ey)=(xe,xe)$, therefore $p_e$
maps $H$ to $\diag(C_G(e)e)$. This mapping is surjective since for $x
\in C_G(e)$, we have $xe=ex$ therefore $(x,x)\in H$ and
$p_e(x,x)=(xe,xe)$. Furthermore, $R_uP(e) \times R_uP^-(e) \subset
\ker(p_e) \subset H$. All this implies our result: let $I$ be such
that $P_I=P(e)$ and $P_I^-=P^-(e)$, let $J\subset I$ be maximal such
that $(L_J,L_J) \subset H(e) = \ker ( L(e) \to C_G(e)e )$ and let $K$
be the complement of $J$ in $I$. The subsets $J$ and $K$ are
orthogonal (since the morphism $L(e) \to C_G(e)e,\ x\mapsto xe$ is the quotient of a finite cover of $L(e)$ by some semi-simple factors and a subgroup of the centre). Furthermore the group $H(e)$ satisfies $(L_J,L_J) \subset H(e)
\subset (L_J,L_J)Z_K$ where $Z_K$ is the center of $L_K$. We have
$H=(U_JH(e) \times U^-_JH(e)) \diag(C_G(e)e)$. Since $C_G(e)e$ is a
quotient by a subgroup contained in $H(e)$ of $L_K$ this concludes
the proof of the first assertion. 

For the second assertion, we use \cite[Theorem 5.18]{renner}: for $s$
a simple reflection, the inclusion $(G \times G)\cdot e \subset
\overline{BsB^-}$ holds if and only if $se = es = e$ \emph{i.e.} if
and only if $(s,1)$ and $(1,s)$ are in $H$. This happens if and only
if $s \in J$. 
\end{proof}

Let $\pi : \Xt \to X$ be a morphism of $G$-embeddings. Let $\xt \in
\Xt$ and $x = \pi( \xt )$. Let $\Ot =(G \times G)\cdot \xt$ and
$\Omega = (G \times G) \cdot x$. We denote by $\Ht$ and $H$ the
stabiliser of $\xt$ and $x$ respectively. There is an inclusion $\Ht
\subset H$. Let $\It,\Jt,\Kt$ and $I,J,K$ be the subsets of $\Delta$
corresponding to $\Ht$ and $H$ according to the previous proposition.  

\begin{lemma}
\label{lemme-H-Ht}
Let $H$ and $\Ht$ be as above. 

1. The groups $\Ht$ and $H$ are simultaneously conjugate to $(U_\It
H_\Jt \times U^-_\It H_\Jt) \diag L_\Kt$ and $(U_I H_J \times U^-_I
H_J) \diag L_K$ with $(L_\Jt,L_\Jt) \subset H_\Jt \subset
(L_\Jt,L_\Jt)Z_\Kt$ and $(L_J,L_J) \subset H_J \subset (L_J,L_J)Z_K$. 

2. We have the inclusions $K \subset \Kt$, $\Jt \subset J$ and $\It
\subset I$. 

3. Assume that $\Xt$ is toroidal, that $\pi$ is proper and that $\Ot$
is closed in $\pi^{-1}(\Omega)$. Then $\It = \Kt = K$ and $\Jt =
\emptyset$. 
\end{lemma}

\begin{proof}
1. Choose $\xt$ a $\diag(T)$-fixed point in $\Ot$. Then $\xt$ is unique up to $T \times T$ action and the same holds for $x = \pi(\xt)$. The result follows from the former proposition since the stabilisers of $\xt$ and $x$ are of the desired for up to conjugation in $T \times T$.

2. The inclusions $\Jt \subset J$ and $\It \subset I$ follow from the
inclusion $\Ht \subset H$. Let $\a$ be a positive root not contained
in the root system generated $\Kt$. Then $\a$ is either a root of
the root system generated by $\Jt$ or we have the inclusion $U_\a
\subset U_\It$. In the first case $\a$ cannot be a root of the root
system generated by $K$. In the second case, we have the inclusions
$U_\a\times\{1\} \subset \Ht \subset H$. It follows, that $\a$ cannot
be a root of the root system generated by $K$. The inclusion $K
\subset \Kt$ follows. 

(\i\i\i) The fact that $\Jt = \emptyset$ and $\Kt = \It$ follows from
the former proposition. With the above assumption, the map
$\pi^{-1}(\Omega) \to \Omega$ is proper and $\Ot$ is closed in
$p^{-1}(\Omega)$. In particular, the map $\Ot \to \Omega$ is
proper. But according to the above proposition, its fibers are
isomorphic to the contracted product $(L_J \times L_J ) \times ^{P_\It
  \cap L_J \times P_\It^- \cap L_J} L$ where $L$ is a quotient of
$L_{\It \cap J}$ by a central subgroup. Since the fiber
is proper it follows that $\It \cap J = \emptyset$.  
\end{proof}

\section{Generalised Schubert and Richardson varieties}

\subsection{Definition and first properties}

Let $X$ be a $G$-embedding. We describe the $B \times B$-orbits and
the $B^- \times B^-$-orbits in any $G \times G$-orbit. Since the $B
\times B$-orbits and the $B^- \times B^-$-orbits are contained in $G
\times G$-orbits, we may fix such an orbit $\Omega$ and according to
Proposition \ref{prop-stab}, there is an element $h_\Omega \in \Omega$
such that the stabiliser $H$ of $h_\Omega$ is of the form  $(U_I H_J
\times U^-_I H_J) \diag L_K$ where $(L_J,L_J) \subset H_J \subset
(L_J,L_J)Z_K$ and $Z_K$ is the center of $L_K$. 

\begin{defn}
\label{defi-gen-rich}
Let $\Omega$, $H$ and $h = h_\Omega \in \Omega$ be as above. Let
$u,v,w,x\in W$. 

1. We denote by $p_\Omega : \Omega \to G/P_I \times G/P_I^-$ the
morphism induced by the inclusion $H \subset P_I \times P_I^-$. 

2. The generalised Schubert cell $\Xo_{u,v}(\Omega)$ is the $B\times
B$-orbit $(Bu\times Bv)\cdot h$. The generalised Schubert variety
$X_{u,v}(\Omega)$ is the closure of $\Xo_{u,v}(\Omega)$ in $X$. 

3. The generalised opposite Schubert cell $\Xo^{w,x}(\Omega)$ is the
$B^-\times B^-$-orbit given by
$(w_0,w_0)\cdot\Xo_{w_0w,w_0x}(\Omega)=(B^-w\times B^-x)\cdot h$. The
generalised opposite Schubert variety $X^{w,x}(\Omega)$ is the closure
of $\Xo^{w,x}(\Omega)$. 

4. The generalised open Richardson variety $\Xo_{u,v}^{w,x}(\Omega)$
is defined as the intersection $\Xo_{u,v}(\Omega)\cap
\Xo^{w,x}(\Omega)$. The generalised Richardson variety
$X_{u,v}^{w,x}(\Omega)$ is defined as the intersection
$X_{u,v}(\Omega)\cap X^{w,x}(\Omega)$. 
\end{defn}

Let $I$, $J$ and $K$ as above. Recall that $u \in W$ can be written $u
= u^I u_I = u^I u_J u_K$. 

\begin{prop}
\label{prop-proj-et-schub}
Let $\Omega$ be a $G \times G$-orbit and let $h \in \Omega$ be as above.

1. The $B \times B$ orbits in $\Omega$ are the generalised Schubert cells. 

2. We have $\Xo_{u,v}(\Omega)=\Xo_{u',v'}(\Omega)$ if and only if $u^I
= {u'}^I$, $(v^{I})^{-1}=({v'}^I)^{-1}$ and
$u_K(v_K)^{-1}=u'_K(v'_K)^{-1}$.  

3. We have $p_\Omega(\Xo_{u,v}(\Omega)) = BuP_I/P_I \times
BvP_I^-/P_I^-$. 

4. The $B^- \times B^-$ orbits in $\Omega$ are the opposite
generalised Schubert cells.  

5. We have $\Xo^{w,x}(\Omega)=\Xo^{w',x'}(\Omega)$ if and only if $w^I
= {w'}^I$, $(x^{I})^{-1}=({x'}^I)^{-1}$ and
$w_K(x_K)^{-1}=w'_K(x'_K)^{-1}$.  

6. We have $p_\Omega(\Xo^{w,x}(\Omega)) = B^- wP_I/P_I \times B^- x
P_I^-/P_I^-$. 
\end{prop}

\begin{proof}
This follows from \cite[Lemma 1.2]{brion-behabiour} since the orbit
$\Omega$ is induced from a quotient $L'$ of $L_K$ in the following
sense $\Omega \simeq (G \times G)^{P_I \times P_I^-}L'$. 
\end{proof}

\begin{example}
\label{exam-1}
In general $X_{u,v}^{w,x}(\Omega)$ is neither irreducible nor
equidimensional. We will however prove in Propposition \ref{prop-CM}
that for $X$ smooth and toroidal, the variety $X_{u,v}^{w,x}(\Omega)$
is irreducible. 

Let $X$ be $\bP(M_3(k))$ where $M_3(k)$ is the vector space of $3
\times 3$ matrices. The group $G \times G$ with $G = {\rm PGL}_3(k)$
acts on $X$ by $(P,Q) \cdot A = PAQ^{-1}$. Let $B$ be the image of the
subgroup of upper-triangular matrices in $G$ and $B^-$ be the image in
$G$ of lower-triangular matrices. For $A \in M_3(k)$, denote by $C_1$,
$C_2$ and $C_3$ the columns of $A$. 

The $G \times G$-orbits are indexed by the rank. Let $\Omega_2$ be the
orbit of matrices of rank $2$. Denote by $X_{1,2}$ and $X_{2,3}$ the
closed subsets given by the equations $C_1 \wedge C_2 = 0$ and $C_2
\wedge C_3 = 0$. The intersections $\overline{\Omega_2} \cap X_{1,2}$ and
$\overline{\Omega_2} \cap X_{2,3}$ are easily seen to be $B \times B$
and $B^- \times B^-$-generalised Schubert varieties. Denote them by
$X_{u,v}(\Omega_2)$ and $X^{w,x}(\Omega_2)$. Let
$X_{u,v}^{w,x}(\Omega_2)$ be the corresponding generalised Richardson
variety. One easily checks that
$X_{u,v}^{w,x}(\Omega_2)$ is the union   
$$X_{u,v}^{w,x}(\Omega_2) = \{ [A] \in X \ | \ C_2 = 0\} \cup \{ [A]
\in X \ | \ {\rm rk}(A) = 1\}.$$ 
This is the decomposition of $X_{u,v}^{w,x}(\Omega_2)$ in irreducible
components. The dimensions of these components are $5$ and
$4$. Therefore $X_{u,v}^{w,x}(\Omega_2)$ is neither irreducible nor
equidimensional.
\end{example}

\begin{prop}
\label{prop-open-smooth}
Let $\Omega$ and $h \in \Omega$ be as above. Let $u,v,w,x \in W$. The
variety $\Xo_{u,v}^{w,x}(\Omega)$ is irreducible and smooth.  
\end{prop}

\begin{proof}
We follow the proof of the same result for rational projective
homogeneous spaces. Let $I,J,K$ such that $H = \textrm{Stab}(h) =
(U_IH_J\times U_I^-H_J)\diag(L_K)$ with $(L_J,L_J) \subset H_J \subset
(L_J,L_J)Z_K$. There is an open dense subset of $\Omega$ given by
$(B^-\times B)\cdot h$. We translate this subset in a neighborhood 
$(wB^-w^{-1}w\times xBx^{-1}x)\cdot h$ of $(w,x)\cdot h$. This
neighborhood contains the $B^-\times B^-$-orbit $(B^-w \times B^-x)\cdot
h_K$. In what follows, we write, for $E$ a subset of $G$ and $\a$ a
root of $(G,T)$: $\a \in E$ for $U_\a \subset E$. We have an
isomorphism given by the action  
$$U_{w,x}\times
(B^-w \times B^-x)\cdot h \simeq   (wB^-w^{-1}w\times xBx^{-1}x)\cdot h$$
with
$$U_{w,x} = \prod_{\a > 0,\ w^{-1}(\a) \in U_I^- \cup B_K^-} U_\a \times
\prod_{\beta > 0,\ x^{-1}(\beta) \in U_I \cup B_K}U_\beta.$$
Intersecting with $\Xo_{u,v}(\Omega)$ which is stable under $U_{w,x}$ we get
$$U_{w,x}\times
\Xo^{w,x}_{u,v}(K) \simeq (wB^-w^{-1}w\times xBx^{-1}x)\cdot h \cap \Xo_{u,v}(K).$$
Since $\Xo_{u,v}(K)$ is irreducible and smooth, the same holds for the right hand side (which is an open subset of $\Xo_{u,v}(K)$) and therefore $\Xo^{w,x}_{u,v}(K)$ is irreducible and smooth.
\end{proof}

\begin{lemma}
The intersection $X_{u,v}^{w,x}(\Omega) \cap \Omega$ is the 
closure of the cell $\Xo_{u,v}^{w,x}(\Omega)$ and is irreducible.
\end{lemma}

\begin{proof}
The variety $X_{u,v}^{w,x}(\Omega) \cap \Omega$ is a union of
intersections $\Xo_{u',v'}(\Omega)\cap\Xo^{w',x'}(\Omega)$  where
$\Xo_{u',v'}(\Omega)$ are the generalised Schubert cells contained in 
$X_{u,v}(\Omega) \cap \Omega$ and $\Xo^{w',x'}(\Omega)$ are the
generalised opposite Schubert cells contained in $X^{w,x}(\Omega) \cap
\Omega$. 

In the orbit $\Omega$, since these Schubert cells are stable for
opposite Borel subgroups of $G\times G$, they are in general position
and thefore intersect properly (see \cite{kleiman:transversality}). In
particular $X_{u,v}^{w,x}(\Omega) \cap \Omega$ contains a unique intersection
$\Xo_{u',v'}(\Omega)\cap\Xo^{u',v'}(\Omega)$ of codimension
$\codim_{{\Omega}} X_{u,v}(\Omega) +
\codim_{{\Omega}} X^{w,x}(\Omega)$: the generalised open
Richardson variety $\Xo_{u,v}^{u,v}(\Omega)$. Since $\Omega$ is
smooth, it follows from \cite[Lemma page 108]{fulton-pragacz} that the
codimension of any irreducible component of $X_{u,v}^{w,x}(\Omega) \cap \Omega$
in ${\Omega}$ is at least $\codim_{{\Omega}}
X_{u,v}(\Omega) + \codim_{{\Omega}} X^{w,x}(\Omega)$. Thus
$X_{u,v}^{w,x}(\Omega) \cap \Omega$ is the  closure of
$\Xo_{u,v}^{w,x}(\Omega)$ and is irreducible.
\end{proof}

\begin{lemma}
\label{lemma-proj-et-schub}
Let $\Omega$ and $h \in \Omega$ be as above. Let $u,v,w,x \in W$. The
closure of the image $p_\Omega(\Xo_{u,v}^{w,x}(\Omega))$ is a product
of projected Richardson varieties in $G/P_I \times G/P_I^-$. 
\end{lemma}

\begin{proof}
We shall see in Section \ref{sec-frob} that all the generalised
Schubert cells, varieties, opposite cells and opposite varieties are
$B \times B$-canonically split for the same splitting. In follows that
all the generalised (open) Richardson varieties are also $B \times
B$-canonically split and the closure of their images
$p_\Omega(\Xo_{u,v}^{w,x}(\Omega))$ are again $B \times
B$-canonically split. Applying \cite[Theorem 5.1]{KLS}, these
varieties are products of projected Richardson varieties.
\end{proof}

\begin{example}
In general $p_\Omega(\Xo_{u,v}^{w,x}(\Omega))$ is not a product of
Richardson variety or even the intersection of opposite $B \times B$
and $B^- \times B^-$-orbits. We will however prove in Proposition
\ref{prop-proj-et-rich} that for $X$ toroidal, the variety
$p_\Omega(X_{u,v}^{w,x}(\Omega))$ is a product of Richardson variety. 

Let $X$ be $\bP(M_4(k))$ where $M_4(k)$ is the vector space of
$4 \times 4$ matrices. The group $G \times G$ with $G = {\rm
  PGL}_4(k)$ acts on $X$ by $(P,Q) \cdot A = PAQ^{-1}$. Let $B$ be the
image of the subgroup of upper-triangular matrices in $G$ and $B^-$ be
the image in $G$ of lower-triangular matrices. For $A \in M_3(k)$,
denote by $C_1$, $C_2$, $C_3$ and $C_4$ the columns of $A$. Let
$(e_1,e_2,e_3,e_4)$ be the canonical basis of $k^4$. 

The $G \times G$-orbits are indexed by the rank. Let $\Omega_2$ be the
orbit of matrices of rank $2$. We have the structure map $p_{\Omega_2}
: \Omega_2 \to \G(2,4) \times \G(2,4)$ defined by $p_{\Omega_2}([A]) =
( \ker A , \im A)$. Here $\G(2,4)$ denotes the Gra{\ss}mann variety of
lines in $\bP^3$. The fiber $p_{\Omega_2}^{-1}(V_2,W_2)$ is the open
subset of $\bP\Hom(k^4/V_2,W_2)$ of invertible elements. Let $\Theta$
be the dense $B \times B$-orbit in $\G(2,4) \times \G(2,4)$ and let
$\Theta^0$ be the dense $B^- \times B^-$-orbit. One easily checks that 
$\{[A] \in \Omega_2\ | \ p_{\Omega_2}([A]) \in \Theta \textrm{ and } 0
\neq A(e_1) \in \scal{e_1,e_2,e_3} \}$  and $\{ [A] \in \Omega_2\ |
\ p_{\Omega_2}([A]) \in \Theta^0 \textrm{ and } 0 \neq A(e_4) \in
\scal{e_2,e_3,e_4} \}$ are irreducible and $B \times B$ resp. $B^-
\times B^-$-stable and therefore contain dense $B \times B$ and $B^-
\times B^-$-orbits that we denote by $\Xo_{u,v}(\Omega_2)$ and
$\Xo^{w,x}(\Omega_2)$.  

We claim that $p_{\Omega_2}(\Xo_{u,v}^{w,x}(\Omega_2))$ is dense in
but different from $\Theta \cap \Theta^0$.  

Let $(V_2,W_2) \in \Theta \cap \Theta^0$ such that $V_2 \cap
\scal{e_1,e_4} = 0$ and $W_2 \cap \scal{e_2,e_3} = 0$. Then the
classes $(\bar e_1, \bar e_4)$ of $e_1$ and $e_4$ in $k^4 / V_2$ form
a basis. Furthermore $W_2 \cap \scal{e_1,e_2,e_3}$ and $W_2 \cap
\scal{e_2,e_3,e_4}$ are in direct sum. Therefore, there is an
isomorphism $[f] \in \bP\Hom(k^4/V_2,W_2)$ with $f(\bar e_1) \in W_2
\cap \scal{e_1,e_2,e_3}$ and $f(\bar e_4) \in W_2 \cap
\scal{e_2,e_3,e_4}$ thus $(V_2,W_2) \in
p_{\Omega_2}(\Xo_{u,v}^{w,x}(\Omega_2))$.  

Let $(V_2,W_2) = (\scal{e_1 + e_4 , e_2 + e_3}) , \scal{e_1 + e_3 ,
  e_2 + e_4}) \in \Theta \cap \Theta^0$. An element $[A] \in
X_{u,v}^{w,x}(\Omega_2)$ with $p_{\Omega_2}([A]) = (V_2 , W_2)$ should
satisfy $0 \neq A(e_1) \in \scal{ e_1 + e_3 }$, $0 \neq A(e_4) \in
\scal{ e_2 + e_4}$ and $A(e_1 + e_4) = 0$. This is impossible. 

Note also that taking inverse images in a toroidal $G$-embedding
dominating $X$ we can also construct examples of this kind for toroidal
varieties. 
\end{example}

\subsection{Generalised Richardson varieties in the smooth toroidal case}

\begin{prop}
\label{prop-CM}
Let $X$ be toroidal and smooth. Generalised Richardson varieties are
irreducible and Cohen-Macau\-lay. 
\end{prop}

\begin{proof}
Let $\Omega$ be a $G \times G$-orbit in $X$. The variety
$X_{u,v}^{w,x}(\Omega)$ is a union of intersections
$\Xo_{u',v'}(\Omega')\cap\Xo^{w',x'}(\Omega')$  where $\Omega'$ is an
$G \times G$-orbit contained in $\overline{\Omega}$, where
$\Xo_{u',v'}(\Omega')$ are the generalised Schubert cells contained in
$X_{u,v}(\Omega) \cap \Omega'$ and where $\Xo^{w',x'}(\Omega')$ are the
generalised opposite Schubert cells contained in $X^{w,x}(\Omega) \cap
\Omega'$.  

In the orbit $\Omega'$, since these Schubert cells are stable for
opposite Borel subgroups of $G\times G$, they are in general position
and thefore intersect  properly (see
\cite{kleiman:transversality}). In particular $X_{u,v}^{w,x}(\Omega)$
contains a unique intersection
$\Xo_{u',v'}(\Omega')\cap\Xo^{u',v'}(\Omega')$ of codimension
$\codim_{\overline{\Omega}} X_{u,v}(\Omega) +
\codim_{\overline{\Omega}} X^{w,x}(\Omega)$: the generalised open
Richardson variety $\Xo_{u,v}^{u,v}(\Omega)$. Since $X$ is smooth and
toroidal, the orbit closure $\overline{\Omega}$ is also smooth (see
for example \cite[Proposition 6.2.4]{BK}). It follows from \cite[Lemma
  page 108]{fulton-pragacz} that the codimension of any irreducible
component of $\Xo_{u,v}^{w,x}(\Omega)$ in $\overline{\Omega}$ is at
least $\codim_{\overline{\Omega}} X_{u,v}(\Omega) +
\codim_{\overline{\Omega}} X^{w,x}(\Omega)$. Thus
$X_{u,v}^{w,x}(\Omega)$ is the closure of $\Xo_{u,v}^{w,x}(\Omega)$
and is irreducible. The Cohen-Macaulay property again follows from
\cite[Lemma page 108]{fulton-pragacz}.
\end{proof}

\begin{prop}
\label{prop-normal}
Let $X$ be toroidal and smooth. Generalised Richardson varieties are
normal. 
\end{prop}

\begin{proof}
Generalised Richardson varieties are Cohen-Macaulay by Proposition
\ref{prop-CM}. It remains to prove that they are smooth in codimension
one. But by Proposition \ref{prop-open-smooth} the generalised open
Richardson varieties are smooth therefore the non smooth locus is
contained in smaller generalised Richardson varieties. The divisorial
part of the non smooth locus is therefore contained in one of these
smaller generalised Richardson varieties. But since all generalised
Richardson varieties are Frobenius split for the same splitting (see
Section \ref{sec-frob}), their intersection is reduced and therefore
generically smooth. It follows that generalised Richardson varieties
are smooth in codimension one. 
\end{proof}

\begin{defn}
Let $(Y')_{Y' \in \mathcal{Y}}$ be a finite family of closed
irreducible subvarieties of an irreducible variety $Y$. The family
$\mathcal{Y}$ is called a stratification if $Y \in \mathcal{Y}$ and
for $Y',Y'' \in \mathcal{Y}$, the intersection $Y' \cap Y''$ is the
union of subvarieties in $\mathcal{Y}$. 
\end{defn}

\begin{prop}
\label{prop-start}
Let $X$ be toroidal and smooth. Generalised Richardson varieties form
a stratification of $X$. 
\end{prop}

\begin{proof}
Since $X_{u,v}^{w,x}(\Omega)$ is irreducible, this follows from the
fact that $X$ is the disjoint union of the open generalised Richardson
varieties. 
\end{proof}

\begin{prop}
\label{prop-proj-et-rich}
Let $X$ be toroidal and let $\Omega$ and $h \in \Omega$ be as
above. Let $u,v,w,x \in W$. The closure of the image
$p_\Omega(\Xo_{u,v}^{w,x}(\Omega))$ is a product of Richardson
varieties in $G/P_I \times G/P_I^-$. 
\end{prop}

\begin{proof}
Note that the image $p_\Omega(\Xo_{u,v}^{w,x}(\Omega))$ is contained
in the product of Richardson varieties $(\overline{BuP_I}/P_I \cap
\overline{B^-wP_I}/P_I) \times (\overline{BvP_I^-}/P_I^- \cap
\overline{B^-xP_I^-}/P_I^-)$. Furthermore its closure is a product of
projected Richardson varieties so it is enough to prove that the
projections to $G/P_I$ and $G/P_I^-$ of the closure of
$p_\Omega(\Xo_{u,v}^{w,x}(\Omega))$ contain the above Richardson
varieties 
$(\overline{BuP_I}/P_I \cap \overline{B^-wP_I}/P_I)$ and
$(\overline{BvP_I^-}/P_I^- \cap \overline{B^-xP_I1-}/P_I^-)$.

Let $\Omega'$ be a closed $G \times G$-orbit in the closure of
${\Omega}$. Since $X$ is toroidal, the orbit $\Omega'$ is isomorphic
to $G/B \times G/B^-$ and we have a commutative diagram (see for
example \cite[Section 5.5]{he-thomsen} for the fact that $p_\Ot$
extends to the closure of $\Omega$): 
$$\xymatrix{
 & \Omega' \ar@{=}[r]^-{p_{\Ot'}}\ar@{^(->}[d] & G/B \times G/B^- \ar[d] \\
\Omega \ar@{^(->}[r] & \overline{\Omega} \ar[r]^-{p_{\Ot}} & G/P_I
\times G/P_I^-}$$ 
According to \cite[Proposition 6.3]{he-thomsen}, the $B \times
B$-orbit $B u' B/B \times B v' B^-/B^- $ in $\Omega'$ is contained in
$X_{u,v}(\Omega)$ if and only if there exists $a \in W_I$ with $u'
\leq ua$ and $v' \geq va$. The same argument proves that the $B^-
\times B^-$-orbit $B^- w' B/B \times B^- x' B^-/B^-$ in $\Omega'$ is
contained in $X^{w,x}(\Omega)$ if and only if there exists $b \in W_I$
with $w' \geq wb$ and $x' \leq xb$. Let $\pi : G/B \to G/P_I$ and
$\pi_- : G/B^- \to G/P_I^-$. For $a$ such that $u' = ua$ is of maximal
length in $u W_I$ and for $b$ such that $x' = xb$ is of maximal length
in $x W_I$ we have $\pi^{-1}(\overline{BuP_I/P_I}) =
\overline{Bu'B/B}$ and $\pi_-^{-1}(\overline{B^-xP^-_I/P^-_I}) =
\overline{B^-x'B^-/B^-}$. Let $v' = va$ and $w' = wa$, we have that
$(\overline{BuP_I}/P_I \cap \overline{B^-wP_I}/P_I) \times
(\overline{BvP_I^-}/P_I^- \cap \overline{B^-xP_I-}/P_I^-)$ is equal to
$\pi(\overline{Bu'B}/B \cap \overline{B^-w'B}/B) \times
\pi_-(\overline{Bv'B^-}/B^- \cap \overline{B^-x'B-}/B-)$ which is
contained in closure of the image
$p_\Omega(\Xo_{u,v}^{w,x}(\Omega))$. 
\end{proof}

\section{Generalised projected Richardson varieties}

\subsection{Definition and first properties}

Recall the following general result on $G$-embeddings.

\begin{prop}
\label{prop-tor2}
1. For any $G$-embedding $X$, there exists a smooth toroidal
$G$-embedding $\Xt$ and a $G\times G$-equivariant morphism
$\psi:\Xt\to X$.

2. For any $G$-embedding $X$ and toroidal $G$-embeddings $\Xt$ and
$\Xt'$ with $G\times G$-equivariant morphisms $\psi:\Xt\to X$ and
$\psi':\Xt'\to X$, there exists a smooth toroidal embedding $\Xh$ with
$G\times G$-equivariant morphisms $\varphi:\Xh\to \Xt$ and
$\varphi:\Xh\to \Xt$ such that the following diagram is commutative.
$$\xymatrix{\Xh\ar[r]^-{\varphi}\ar[d]_-{\varphi'} & \Xt\ar[d]^-\psi \\
\Xt'\ar[r]^-{\psi'} & X \\}$$
\end{prop}

\begin{proof}
1. This result is proved in \cite[Theorem 6.2.5]{BK}. 

2. Ths result is classical for spherical varieties in general (see
for example \cite{knop}) without the smoothness condition on $X''$ but
using 1. the result follows for $G$-embeddings. 
\end{proof}

\begin{defn}
\label{defi-projected}
Let $X$ be a $G$-embedding, a projected generalised Richardson variety
is the image of a generalised Richardson variety $\Xt_{u,v}^{w,x}(\Omega)$
under an equivariant morphism $\varphi:\Xt\to X$ with $\Xt$ smooth and
toroidal.
\end{defn}

\begin{lemma}
\label{lemme-tor}
Let $\varphi:Y\to X$ be a $G\times G$-equivariant morphism between
smooth toroidal $G$-embeddings. Let $u,v,w,x \in W$.

1. Let $\Omega$ be a $G \times G$-orbit in $X$ then there exists
a $G \times G$-orbit $\Omega'$ in $Y$ such that
$\varphi(\Yo_{u,v}^{w,x}(\Omega'))={\Xo}_{u,v}^{w,x}(\Omega)$ and
$\varphi(Y_{u,v}^{w,x}(\Omega'))={X}_{u,v}^{w,x}(\Omega)$. 

2. Let $\Omega'$ be a $G \times G$-orbit in $Y$ and $\Omega =
\varphi(\Omega')$. Then
$\varphi(\Yo_{u,v}^{w,x}(\Omega'))={\Xo}_{u,v}^{w,x}(\Omega)$ and
$\varphi(Y_{u,v}^{w,x}(\Omega'))={X}_{u,v}^{w,x}(\Omega)$. 
\end{lemma}

\begin{proof}
It is enough to prove 2. since for any $G \times G$-orbit in $X$, there
exists a $G \times G$-orbit $\Omega'$ in $Y$ such that
$\varphi(\Omega') = \Omega$. Write $\Omega' = (G \times G) \cdot \xt$ and
$\Omega = (G \times G) \cdot x$ and let $\Ht$ and $H$ the stabilisers
of $\xt$ and $x$ in $G \times G$. According to Lemma \ref{lemme-H-Ht}, we have
$H = \Ht Z$ for $Z$ a subgroup of $T \times T$. In particular, the map
$\varphi: \Omega' \to \Omega$ is a quotient by $Z$, maps
$\Yo_{u,v}(\Omega')$ to $\Xo_{u,v}(\Omega)$ and
$\varphi^{-1}(\Xo_{u,v}(\Omega))=\Yo_{u,v}(\Omega')$. The same holds for the
opposite Schubert cells: $\varphi(\Yo^{w,x}(\Omega'))=\Xo^{w,x}(\Omega)$ and
$\varphi^{-1}(\Xo^{w,x}(\Omega))=\Yo^{w,x}(\Omega')$. Taking closures, the same
result holds for generalised Schubert varieties and generalised
opposite Schubert varieties. We deduce that
$\varphi(\Yo_{uv}^{w,x}(\Omega'))=\Xo_{uv}^{w,x}(\Omega)$ and taking closures,
the result follows (recall that for $X$ smooth and toroidal,
the variety $X_{u,v}^{w,x}(\Omega)$ is irreducible).
\end{proof}

\begin{cor}
\label{cor-relevt}
Let $\psi:\Xt \to X$ be a morphism of
$G$-embeddings with $\Xt$ smooth and toroidal. Then any projected
generalised Richardson variety is the projection of a
generalised Richardson variety in $\Xt$.
\end{cor}

\begin{proof}
Let $\psi:\Xt\to X$ and $\psi':\Xh\to X$ be two smooth toroidal
variety dominating $X$. Let ${\Xh}_{u,v}^{w,x}(\Omega')$ be a generalised
Richardson variety in $\Xt$, we prove that
$\psi'({\Xh}_{u,v}^{w,x}(\Omega'))$ is also the projection of a
generalised Richardson variety in $\Xt$. Let $\Xt''$ smooth and
toroidal dominating both
$\Xt$ and $\Xh$ as given in Proposition \ref{prop-tor2}. Then
$\varphi({\varphi'}^{-1}({\Xh}_{u,v}^{w,x}(\Omega')))$ is again a
generalised Richardson variety in $\Xt$ and the result follows.
\end{proof}

\subsection{Parabolic induction}

In this subsection we consider the following situation. Let $\cG$ be a
reductive group, $\cT$ be a maximal torus and $\cP$ be a parabolic
subgroup containing $\cT$. Let $\cU$ be the unipotent radical of
$\cP$. We denote by $\cW,\cW_P$ the Weyl groups of $(\cG,\cT)$ and
$(\cP,\cT)$ and by $\cL$ the Levi subgroup of $\cP$ containing
$\cT$. Let $\cB$ be a Borel subgroup of $\cG$ with $\cT \subset \cB
\subset \cP$ and let $\cB^-$ be the opposite Borel subgroup with
repect to $\cT$. We write $\cW^\cP$ for the set of minimal length
representatives of $\cW/\cW_P$. 

Let $\cH$ be a spherical subgroup of $\cG$ contained in $\cP$ such
that $\cU \subset \cH$, let $\cX = \cG/\cH$ and $\cY = \cG/\cP$. We
have $\cX \simeq \cG \times^\cP \cP/\cH$ and $\cP/\cH \simeq \cL/\cL
\cap \cH$. The quotient $\cP/\cH$ is thus a $\cL$-spherical
variety. Let $p : \cX \to \cY$ be the natural projection. By
\cite[Lemma 1.2]{brion-behabiour}), any $\cB$-orbit of $\cX$ is of the
form $\cB \lambda \cO$ for $\lambda \in \cW^P$ and $\cO$ a $\cB_\cL =
\cB \cap \cL$-orbit in $\cP/\cH$ and any $\cB^-$-orbit of $\cX$ is of
the form $\cB^- \mu \cO^-$ for $\mu \in \cW^P$ and $\cO^-$ a
$\cB^-_\cL = \cB^- \cap \cL$-orbit in $\cP/\cH$. 

\begin{lemma}
Let $b \in \cB$ and $b_- \in \cB^-$, let $\lambda,\mu \in \cW^\cP$,
let $\cO = \cB_\cL \nu \cdot \cH$ with $\nu \in \cL$ a $\cB_\cL$-orbit
in $\cP/\cH$ and let $\cO^- = \cB^-_\cL \xi \cdot \cH$ with $\xi \in
\cL$ a $\cB^-_\cL$-orbit in $\cP/\cH$. 

1. We have $\cB \lambda \cO = \cB \lambda \nu \cdot \cH$ and $\cB^-
\mu \cO^- = \cB^- \mu \xi \cdot \cH$. 

2. The intersection $p^{-1}(b \lambda \nu \cdot \cP) \cap \cB \lambda
\cO$ is $b \lambda \nu (\cP \cap \cB^{(\lambda \nu)^{-1}}) \cdot
\cH$. 

3. The intersection $p^{-1}(b_- \mu \xi \cdot \cP)) \cap \cB^- \mu
\cO^-$ is $b_- \mu \xi (\cP \cap (\cB^{-})^{(\mu \xi)^{-1}}) \cdot
\cH$. 

4. Assume that $b \lambda \nu \cdot \cP = b_- \mu \xi \cdot \cP$. Let
$\zeta = (b \lambda \nu)^{-1}(b_- \mu \xi)$. Then we have $$p^{-1}(b
\lambda \nu \cdot \cP) \cap \cB \lambda \cO\cap \cB^- \mu \cO^- =  b
\lambda \nu (\cP \cap \cB^{(\lambda \nu)^{-1}}) \cdot \cH \cap b
\lambda \nu (\cP \cap \zeta (\cB^{-})^{(\mu \xi)^{-1}}) \cdot \cH.$$ 

5. Under the isomorphism $b \lambda \nu \cP \cdot \cH \simeq \cP/\cH
\simeq \cL/(\cL \cap \cH)$, we have  
$$p^{-1}(b \lambda \nu \cdot \cP) \cap \cB \lambda \cO\cap \cB^- \mu
\cO^- \simeq \left( {\cB_\cL}^{\nu^{-1}} \cdot (\cL \cap \cH) \right)
\cap \left( \zeta_l ({\cB_\cL}^{-})^{\xi^{-1}}) \cdot (\cL \cap \cH)
\right).$$ 
where $\zeta = \zeta_l \zeta_u$ with $\zeta_l \in \cL$ and $\zeta_u
\in \cU$. 
\end{lemma}

\begin{proof}
1. Let $b_1 \lambda b_2 \nu \cdot \cH \in \cB \lambda \cO$ with $b_1
\in \cB$ and $b_2 \in \cB_\cL$. Since $\lambda \in \cW^\cP$, we have
$\lambda b_2 \lambda^{-1} \in \cB$ and $b_1 \lambda b_2 \lambda^{-1}
\lambda \nu \cdot \cH \in \cB \lambda \nu \cdot \cH$. This proves $\cB
\lambda \cO = \cB \lambda \nu \cdot \cH$. A similar argument proves
the second equality. 

2. We have $p^{-1}(b \lambda \nu \cdot \cP) \cap \cB \lambda \nu \cdot
\cH = b \lambda \nu \cP \cdot \cH \cap \cB \lambda \nu \cdot \cH$ and
the equality follows. 

3. A similar argument as in 2. proves the result.

4. Follows from 2., 3. and the equality 
$$b_- \mu \xi (\cP \cap (\cB^{-1})^{(\mu \xi)^{-1}}) \cdot \cH = b
\lambda \nu (\cP \cap \zeta (\cB^{-1})^{(\mu \xi)^{-1}}) \cdot \cH.$$ 

5. Let $p \in \cP$, then there is a unique decomposition $p = p_l p_u$
with $p_l \in \cL$ and $p_u \in \cU$ and the map $\cP / \cH \to \cL /
\cL \cap \cH$ is given by $p \cdot \cH \mapsto p_l \cdot (\cL \cap
\cH)$. Furthermore, the map $p \mapsto p_l$ is multiplicative and maps
$\cH$ to $\cL \cap \cH$.

Since $\lambda , \mu \in \cW^\cP$, we have $\cB^{\lambda^{-1}} \cap
\cL = \cB_\cL$ and $(\cB^{-})^{\mu^{-1}} \cap \cL =
\cB^-_\cL$. Furthermore, since $\nu , \xi \in \cL$, we have   
$$\cB^{(\lambda \nu)^{-1}} \cap \cL = {\cB_\cL}^{\nu^{-1}} \textrm{
  and } (\cB^-)^{(\mu \xi)^{-1}} \cap \cL = {\cB^-_\cL}^{\xi^{-1}}.$$
Now for $p \in \cP \cap \cB^{(\lambda \nu)^{-1}}$ and $q \in \cP
\cap {\cB^-}^{(\mu \xi)^{-1}}$ we have $p_l \in \cB^{(\lambda
  \nu)^{-1}} \cap \cL = {\cB_\cL}^{\nu^{-1}}$ and $q_l \in
     {\cB^-}^{(\mu \xi)^{-1}} \cap \cL = {\cB^-_\cL}^{\xi^{-1}}$.

Let $b \lambda \nu p \cdot \cH \in p^{-1}(b \lambda \nu \cdot \cP)
\cap \cB \lambda \cO\cap \cB^- \mu \cO^-$. Then $b \lambda \nu p \cdot
\cH$ is mapped to $p_l \cdot (\cL \cap \cH)$ in $\cL/(\cL \cap
\cH)$. Furthermore, according to 4., there are elements $h_1,h_2 \in
\cH$ such that $ph_1 \in (\cB)^{(\lambda \nu)^{-1}}$ and $\zeta^{-1} p
  h_2 \in {\cB^-}^{(\mu \xi)^{-1}}$. Then $p_l(h_1)_l = (ph_1)_l \in
  {\cB_\cL}^{\nu^{-1}}$ and $\zeta^{-1}_l p_l (h_2)_l = (ph_2)_l \in
  {\cB_\cL^-}^{\xi^{-1}}$. Since $(h_1)_l,(h_2)_l \in (\cL \cap \cH)$,
  we have $p_l \cdot (\cL \cap \cH) = p_l (h_1)_l \cdot (\cL \cap \cH) \in
  {\cB_\cL}^{\nu^{-1}} \cdot (\cL \cap \cH)$ and $p_l \cdot (\cL \cap
  \cH) = p_l (h_2)_l \cdot (\cL \cap \cH) \in \zeta_l
  ({\cB_\cL}^{-})^{\xi^{-1}} \cdot (\cL \cap \cH)$. Therefore
$$p_l \cdot (\cL \cap \cH) \in {\cB_\cL}^{\nu^{-1}} \cdot (\cL \cap
  \cH) \cap \zeta_l ({\cB_\cL}^{-})^{\xi^{-1}}) \cdot (\cL \cap
  \cH).$$
The converse inclusion is easy.
\end{proof}

We apply the above result to the following situation. Let $X$ be a
$G$-embedding and $\Omega = (G \times G) \cdot x$ such that the
stabiliser $H$ of $x$ is as given in Proposition \ref{prop-stab}: $H =
(U_IH_J \times U_I^- H_J)\diag(L_K)$. Let $\cG = G \times G$, $\cT = T
\times T$, $\cB = B \times B$, $\cH = H$ and $\cP = P_I \times P_I^-$
. We have $\cU = U_I \times U_I^- \subset \cH$ and $\cL = L_I \times
L_I = (L_J \times L_K) \times (L_J \times L_K)$. 

\begin{cor}
1. The fibers of the map $p_\Omega:\Omega \to G/P_I \times G/P_I^-$ are
isomorphic to a quotient $L$ of $L_K$ by a central subgroup.

2. The fibers of the restriction $p_\Omega: \Xo_{u,v}^{w,x}(\Omega)
\to p_\Omega(\Xo_{u,v}^{w,x}(\Omega))$ are isomorphic to the
intersection of $(B_K \times B_K)^{(u_K,v_K)^{-1}} \cdot 1_L$ and a
translate $\zeta (B_K^- \times B_K^-)^{(w_K,x_K)^{-1}} \cdot 1_L$ for some
$\zeta \in L$ (depending on the fiber). 
\end{cor}

\begin{proof}
1. The fibers are isomorphic to $\cP/\cH \simeq \cL/(\cL \cap \cH) =
(L_J \times L_K) \times (L_J \times L_K) / (H_J \times H_J)
\diag(L_K)$. The last term is isomorphic to a quotient of $L_K$ by a
central subgroup (contained in $H_J$).

2. The $B \times B$-orbit $\Xo_{u,v}(\Omega)$ is of the form $\cB
(u,v) \cdot \cH = \cB (u^I,v^I) \cB_L (u_I,v_I) \cdot \cH$. A similar
statement holds for $\Xo^{w,x}(\Omega)$. Applying the above Lemma,
we get that the fibers are isomorphic to the
intersection of $(B_I \times B_I)^{(u_I,v_I)^-1} \cdot 1_L$ and a
translate $\zeta (B_I^- \times B_I)^{(w_I,x_I)^-1} \cdot 1_L$ for some
$\zeta \in L$. Since $L$ is a quotient of $L_J \times L_K$ by $(H_J
\times H_J)\diag(L_K)$ the result follows.
\end{proof}

\begin{remark}
Note that in $L$ we have $(B_K \times B_K)^{(u_K,v_K)^{-1}} \cdot 1_L
= B_K^{u_K^{-1}}B_K^{v_K^{-1}}$ and $\zeta (B_K^- \times
B_K^-)^{(w_K,x_K)^{-1}} \cdot 1_L = \zeta \cdot
({B_K^-}^{w_K^{-1}}{B_K^-}^{x_K^{-1}})$. 
\end{remark}

Let $\varphi : \Xt \to X$ be a morphism of $G$-embeddings with $\Xt$
toroidal. According to Lemma \ref{lemme-H-Ht} and Proposition
\ref{prop-stab}, there exists $\xt \in \Xt$ and $x = \varphi(\xt)$
such that if $\Ht$ and $H$ are the stabilisers of $\xt$ and $x$, then  
$$\Ht = (U_\It H_\Jt \times U^-_\It H_\Jt) \diag L_\Kt \textrm{ and }
H = (U_I H_J \times U^-_I H_J) \diag L_K$$  
with $(L_\Jt,L_\Jt) \subset H_\Jt \subset (L_\Jt,L_\Jt)Z_\Kt$,
$(L_J,L_J) \subset H_J \subset (L_J,L_J)Z_K$ and $K \subset \Kt = \It
\subset I$. Note that $\It = K \cup (\It \cap J)$ and that $K$ and
$\It \cap J$ are orthogonal Applying the above
result we get. 

\begin{cor} 
\label{cor-fibres}
Let $\Ot = (G \times G) \cdot \xt$ and $\Omega = (G \times G) \cdot x$.

1. There is a commutative diagram
$$\xymatrix{\Ot \ar[r]^-\varphi \ar[d]_-{p_\Ot} & \Omega \ar[d]^-{p_\Omega} \\
G/P_\It \times G/P_\It^- \ar[r] & G/P_I \times G/P_I.}$$
The fibers of $p_\Ot$ and $p_\Omega$ are isomorphic to quotients of
$L_K \times L_{\It \cap J}$ and $L_K$ by central subgroups. The
morphism between these fibers induced by $\varphi$ is the morphism
induced by the the first projection. 

2. Let $u,v,w,x \in W$. There is a commutative diagram
$$\xymatrix{{\Xto}_{u,v}\ \!\!\!\!\!\!\!\!\!^{w,x}(\Ot)
  \ar[r]^-\varphi \ar[d]_-{p_\Ot} & \Xo_{u,v}^{w,x}(\Omega)
  \ar[d]^-{p_\Omega} \\ 
p_\Ot({\Xto}_{u,v}\ \!\!\!\!\!\!\!\!\!^{w,x}(\Ot)) \ar[r] &
p_\Omega(\Xo_{u,v}^{w,x}(\Omega)),}$$ 
with vertical fibers isomorphic to
$$\left( B_K^{u_K^{-1}}B_K^{v_K^{-1}} \cap \widetilde{\zeta}_K \cdot
{B_K^-}^{w_K^{-1}}{B_K^{-}}^{x_K^{-1}} \right) \times \left( B_{\It
  \cap J}^{u_{\It \cap J}^{-1}}B_{\It \cap J}^{v_{\It \cap J}^{-1}}
\cap \widetilde{\zeta}_{\It \cap J} \cdot {B_{\It \cap J}^-}^{w_{\It
    \cap J}^{-1}}{B_{\It \cap J}^{-}}^{x_{\It \cap J}^{-1}} \right)$$ 
and 
$$B_K^{u_K^{-1}}B_K^{v_K^{-1}} \cap {\zeta}_K \cdot
{B_K^-}^{w_K^{-1}}{B_K^{-}}^{x_K^{-1}}.$$ 
Furthermore, the morphism between these fibers induced by $\varphi$ is
the morphism induced by the first projection.
\end{cor}

\subsection{Stratification}

Let $X$ be a proper $G$-embedding. In this subsection, we prove that
the projected generalised Richardson varieties in $X$ form a
stratification. For this, according to Corollary \ref{cor-relevt},  we
can fix a smooth toroidal variety $\Xt$ together with a proper $G
\times G$-equivariant morphism $\varphi : \Xt \to X$. All the
projected generalised Richardson varieties are of the form
$\varphi(\Xt_{u,v}^{w,x}(\Ot))$ for some orbit $\Ot$ and elements
$u,v,w,x \in W$. 

\begin{defn}
1. For each $G \times G$-orbit in $X$, we choose a $G \times G$-orbit
$\Ot$ in $\Xt$ such that $\Ot$ is minimal in
$\varphi^{-1}(\Omega)$. We define $I,J,K,\It,\Jt,\Kt$ as the subsets
of simple roots such that  
\begin{itemize}
\item $\Omega \simeq (G\times G) / H$ with $H = (U_I H_J \times U^-_I
  H_J) \diag L_K$ and $(L_J,L_J) \subset H_J \subset (L_J,L_J)Z_K$. 
\item $\Ot \simeq (G\times G) / \Ht$ with $\Ht = (U_\It H_\Jt \times
  U^-_\It H_\Jt) \diag L_\It$ and $H_\Jt \subset Z_\It$. 
\end{itemize}
Recall from Lemma \ref{lemme-H-Ht} that we have $\Jt = \emptyset$ and
$K = \Kt = \It$ and that the roots in $J$ and $K$ are orthogonal. We
write $\psi : G/P_K \times G/P_K^- \to G/P_I \times G/P_I^-$. 
  
2. The set $\mathfrak{R}_X$ is the set of tuples $(\Omega,u,v,w,x)$
with $\Omega$ a $G \times G$-orbit of $X$ and $u,v,w,x \in W$ with $u
= u^I$ and $x = x^I$. 

3. For $\Omega$ a $G \times G$-orbit in $X$ and $u,v,w,x \in W$, we
define 
\begin{itemize}
\item $\Rto_{u,v}\!\!\!\!\!\!\!^{w,x}(\Ot) = (BuP_K \cap B^-wP_K)/P_K
  \times (BvP_K^- \cap B^-xP_K^-)/P_K^-$ and $\Rt_{u,v}^{w,x}(\Ot)$
  its closure. 
\item $\Ro_{u,v}^{w,x}(\Omega) = (BuP_I \cap B^-wP_I)/P_I \times
  (BvP_I^- \cap B^-xP_I^-)/P_I^-$ and $R_{u,v}^{w,x}(\Omega)$ its
  closure. 
\item $\PRo_{u,v}^{w,x}(\Omega) =
  \psi(\Rto_{u,v}\!\!\!\!\!\!\!^{w,x}(\Ot))$ and
  $PR_{u,v}^{w,x}(\Omega) = \psi(\Rt_{u,v}^{w,x}(\Ot))$ its closure. 
\item $\Pi_{u,v}^{w,x}(\Omega) = \varphi(\Xt_{u,v}^{w,x}(\Ot))$ and
  $\Pio_{u,v}^{w,x}(\Omega) =
  \varphi(\Xto_{u,v}\!\!\!\!\!\!\!^{w,x}(\Ot))$. 
\end{itemize}
\end{defn}

\begin{lemma}
\label{lemm-fibre-iso}
Let $(\Omega,u,v,w,x) \in \mathfrak{R}_X$ and $\Ot$ as above.

1. In the commutative diagram
$$\xymatrix{\Ot \ar[r]^-\varphi \ar[d]_-{p_\Ot} & \Omega \ar[d]^-{p_\Omega} \\
G/P_\It \times G/P_\It^- \ar[r]^-\psi & G/P_I \times G/P_I,}$$
for $\yt \in G/P_\It \times G/P_\It^-$, the map $p_\Ot^{-1}(\yt) \to
p_\Omega^{-1}(\psi(\yt))$ induced by $\varphi$ is an isomorphism. 

2. In the commutative diagram
$$\xymatrix{{\Xto}_{u,v}\ \!\!\!\!\!\!\!\!\!^{w,x}(\Ot)
  \ar[r]^-\varphi \ar[d]_-{p_\Ot} & \Xo_{u,v}^{w,x}(\Omega)
  \ar[d]^-{p_\Omega} \\ 
p_\Ot({\Xto}_{u,v}\ \!\!\!\!\!\!\!\!\!^{w,x}(\Ot)) \ar[r] & p_\Omega(\Xo_{u,v}^{w,x}(\Omega)),}$$
for $\yt \in {\Rto}_{u,v}\ \!\!\!\!\!\!\!\!\!^{w,x}(\Ot)$, the map
$p_\Ot^{-1}(\yt) \to p_\Omega^{-1}(\psi(\yt))$ induced by $\varphi$ is
an isomorphism. 
\end{lemma}

\begin{proof}
1. Since $\varphi$ is surjective, the map $p_\Ot^{-1}(\yt) \to
p_\Omega^{-1}(\psi(\yt))$ is surjective. According to Proposition
\ref{prop-stab} and Lemma \ref{lemme-H-Ht}, we can write $\Ot = G
\times G /\Ht$ and $\Omega = G \times G / H$ such that $\Ht = (U_K
\Ht_J \times U_K^- \times \Ht_J) \diag(L_K)$ and $H = (U_I H_J \times
U_I^- \times H_J) \diag(L_K)$ with $\Ht_J \subset Z_K$, $(L_J,L_J)
\subset H_J \subset (L_J,L_J) Z_K$, $I = J \cup K$ and $J$ orthogonal
to $K$. The fibers of $p_\Ot$ and $p_\Omega$ are therefore isomorphic
to $L_K / \Ht_J$ and $L_K / H'$ for $H'$ some subgroup of $Z_K$. It
follows that the map $p_\Ot^{-1}(\yt) \to p_\Omega^{-1}(\psi(\yt))$
induced by $\varphi$ is surjective with fiber isomorphic to the
subgroup $H'/\Ht_J$ of $Z_K$. It also follows that for
$\Xto_{u',v'}(\Ot)$ a $B \times B$-orbit in $\Ot$, the fiber of the
map $\Xto_{u',v'}(\Ot) \to \Xo_{u',v'}(\Omega)$ contains
$H'/\Ht_J$. We prove that this subgroup must be trivial. 

Recall from \cite[Corollary 3.3]{knop-cmh} that if a homogeneous
spherical variety is such that the stabiliser in a Borel subgroup of a
general point is connected, then so is the stabiliser in a Borel
subgroup of any point. In particular, this holds for $G$-embeddings
and their $G \times G$-orbits. Let $x \in \Xo_{u,v}(\Omega)$ let $(B
\times B)_x$ be its stabiliser in $B \times B$. It is connected. Since
it is solvable, it therefore acts with a fixed point $\xt$ on the
fiber $\varphi^{-1}(x) \cap \Ot$ which is closed. The $B \times
B$-orbit of $\xt$, which is of the form $\Xto_{u',v'}(\Ot)$, is
therefore isomorphic to $\Xo_{u,v}(\Omega)$ via $\varphi$. In
particular $H'/\Ht_J$ is trivial. 

2. Follows from 1. and Corollary \ref{cor-fibres}.
\end{proof}

We will need the following result generalising Theorem 3.6 in
\cite{KLS} (see also Theorem 7.1 in \cite{rietsch}).

\begin{lemma}
Let $Q \subset P \subset G$ be parabolic subgroups containing $B$ and
let $p_{Q,P} : G/Q \to G/P$ be the projection. If $\Ro_u^w(Q)$ and
$\Ro_u^w(P)$ denote the open Richardson variety $BuQ/Q \cap B^-wQ/Q$
and $BuP/P \cap B^-wP/P$, then for $u \in
W^P$ and $w \in W$, we have the equality
$$\Ro_u^w(P) = \coprod_{\substack{w' \in W^Q \\ (w')^P = w^P}}
p_{Q,P}(\Ro_u^{w'}(Q)).$$
\end{lemma}

\begin{proof}
Since $p_{Q,P}(R_{u'}^{w'}(Q)) \subset R_{u'}^{w'}(P)$ and
$R_{u'}^{w'}(P) = R_{u'}^{w'}(P)$ for ${u'}^P = u^P$ and ${w'}^P =
w^P$, we have the inclusion:  
$$\bigcup_{\substack{w' \in W^Q \\ (w')^P = w^P}} p_{Q,P}(\Ro_u^{w'}(Q))
\subset \Ro_u^w(P).$$
Consider the commutative diagram
$$\xymatrix{G/B \ar[r]^-{p_{B,Q}} \ar[rd]_-{p_{B,P}} & G/Q
  \ar[d]^-{p_{Q,P}} \\
 & G/P\\}$$
and denote by $\Ro_u^w(B)$ the open Richardson variety $BuB/B \cap
B^-wB/B$ in $G/B$. The same argument as above together with
\cite[Theorem 3.6]{KLS} and the fact that $u = u^P = u^Q$ gives 
$$\Ro_u^{w'}(Q) = \coprod_{\substack{w'' \in W \\ {w''}^Q = {w'}^Q}}
p_{B,Q}(\Ro_u^{w''}(B)) \textrm{ and } \Ro_u^{w}(P) =
\coprod_{\substack{w''' \in W \\ {w'''}^P = w^P}} p_{B,P}(\Ro_u^{w'''}(B)).$$   
Note also that $p_{B,P}(\Ro_u^{w''}(B)) =
p_{Q,P}(p_{B,Q}(\Ro_u^{w''}(B)))$ and that this locally closed
subvarieties are disjoint for $u \in W^P$ fixed (see \cite[Theorem
  3.6]{KLS} again). This implies 
$$p_{Q,P}(\Ro_u^{w'}(Q)) = \coprod_{\substack{w'' \in W \\ (w'')^Q = {w'}^Q}}
p_{B,P}(\Ro_u^{w''}(B)).$$ 
We get
$$\Ro_u^{w}(P) = \coprod_{\substack{w''' \in W \\ {w'''}^P = w^P}}
p_{B,P}(\Ro_u^{w'''}(B)) = \coprod_{\substack{w' \in W^Q \\ {w'}^P =
    w^P}} \coprod_{\substack{w''' \in W \\ {w'''}^Q = {w'}^Q}} 
p_{B,P}(\Ro_u^{w'''}(B))$$    
and the result follows.
\end{proof}

\begin{prop}
\label{prop-strat-gen}
The family $(\Pi_{u,v}^{w,x}(\Omega))_{(\Omega,u,v,w,x) \in
  \mathfrak{R}_X}$ is a stratification of $X$. 
\end{prop}

\begin{proof}
We prove the equality
$$X = \coprod_{(\Omega,u,v,w,x) \in \mathfrak{R}_X}
\Pio_{u,v}^{w,x}(\Omega).$$ 
Let $x \in X$ and let $\Omega$ continaing $x$. Fix a $G \times
G$-orbit $\Ot$ minimal in $\varphi^{-1}(\Omega)$.  
There exist $u,v,w,x \in W$ such that 
$(\Omega,u,v,w,x) \in \mathfrak{R}_X$ and $x \in
\Xo_{u,v}^{w,x}(\Omega)$. Let $y = p_\Omega(x)$. We have $y \in
\Ro_{u,v}^{w,x}(\Omega)$ and by the former lemma there exist uniquely
determined elements $v',w' \in W^{K}$ with ${v'}^I = v^I$ and
${w'}^I = w^I$ such that $y \in \psi(\Rto_{u,v'}^{w',x}(\Ot))$. Let
$\yt \in \Rto_{u,v'}^{w',x}(\Ot)$ with $\psi(\yt) = y$. Note that
we also have $\Ro_{u,v}^{w,x}(\Omega) =
\Ro_{u,v'}^{w',x}(\Omega)$. Let $v'' = v' v_K$ and $w'' = w' w_K$. We
have $\Xo_{u,v}^{w,x}(\Omega) = \Xo_{u,v''}^{w'',x}(\Omega)$ and
$\Rto_{u,v'}^{w',x}(\Ot) = \Rto_{u,v''}^{w'',x}(\Ot)$. By
Lemma \ref{lemm-fibre-iso}, there exists an element $\xt 
\in \Xto_{u,v''}^{w'',x}(\Ot)$ with $p_\Ot(\xt) = \yt$ and
$\varphi(\xt) = x$. It follows that $x \in
\Pio_{u,v''}^{w'',x}(\Omega)$ with $(\Omega,u,v'',w'',x) \in
\mathfrak{R}_X$ uniquely determined.
\end{proof}

\begin{prop}
\label{prop-smooth-gen}
For $(\Omega,u,v,w,x) \in \mathfrak{R}_X$, the variety
$\Pio_{u,v}^{w,x}(\Omega)$ is smooth. 
\end{prop}

\begin{proof}
Since $u = u^I$ and $x = x^I$, in the commutative diagram
$$\xymatrix{{\Xto}_{u,v}\ \!\!\!\!\!\!\!\!\!^{w,x}(\Ot)
  \ar[r]^-\varphi \ar[d]_-{p_\Ot} & \Xo_{u,v}^{w,x}(\Omega)
  \ar[d]^-{p_\Omega} \\ 
{\Rto}_{u,v}\ \!\!\!\!\!\!\!\!\!^{w,x}(\Ot) \ar[r]^-\psi &
\PRo_{u,v}^{w,x}(\Omega),}$$ 
the map $\psi$ is an isomorphism. By Lemma \ref{lemm-fibre-iso} so is
the map $\varphi$ on its image $\Pio_{u,v}^{w,x}(\Omega)$. Since
$\Xto_{u,v}^{w,x}(\Ot)$ is smooth the result follows. 
\end{proof}

\section{Frobenius splittings}
\label{sec-frob}

\subsection{Existence of a splitting} Let $X$ be a
$G$-embedding, then $X$ admits a $B\times B$-canonical splitting (see
\cite[Theorem 6.2.7]{BK}). In \cite{he-thomsen}, He and Thomsen
exhibit many compatibly split subvarieties of a particular
splitting. We recall their results. Write $D_\alpha$ for
$X_{w_os_\alpha,1}(G)$ and $\Dt_\alpha$ for $X^{w_os_\alpha,1}(G)$
(recall that $G$ is the dense orbit in $X$).  

\begin{prop}
\label{thm-HT}
There exists a splitting of $X$ compatibly splitting the irreducible
$G\times G$-divisors, the divisors $(D_\alpha)_{\alpha\in I}$ and the
divisors $(\Dt_\alpha)_{\alpha\in I}$. For $X$ complete, this
splitting is unique. 

This splitting is a $(p-1)$-th power of a global section of
$\omega_X^{-1}$. It compatibly splits the projected generalised
Richardson varieties.
\end{prop}

\begin{proof}
We start with $X$ toroidal. The existence of this splitting (and the
fact that it is a $(p-1)$-th power of a global section of
$\omega_X^{-1}$) is proved in \cite[Theorem 6.2.7]{BK}. The unicity
follows from general arguments: let $\phi$ be a
Frobenius splitting compatibly splitting the irreducible $G\times
G$-divisors, the divisors $(D_\alpha)_{\alpha\in I}$ and the divisors
$(\Dt_\alpha)_{\alpha\in I}$. This splitting is given by a global
section $\sigma$ of $\omega_X^{1-p}$. From \cite[Theorem 1.4.10]{BK}
it follows that $\sigma$ is a global section of 
$$\mathcal{L}=\omega_X^{1-p}\left(-\sum_j(p-1)X(j)-\sum_{\alpha\in
  I}(p-1)(D_\alpha+\Dt_\alpha)\right),$$ 
where the $X(j)$ are the irreducible $G\times G$-divisors on $X$. By
\cite[Proposition 6.2.6]{BK} we have $\mathcal{L}\simeq\cO_X$. The
uniqueness follows.

The second part follows from He and Thomsen's results in
\cite{he-thomsen}. By \cite[Proposition 6.5]{he-thomsen}, all
generalised Schubert varieties and opposite 
Schubert varieties are compatibly split as irreducible components of
intersections of the compatibly split generalised Schubert
divisors. We conclude that all generalised Richardson varieties are
compatibly split. 

By projection, using \cite[Lemma 1.1.8]{BK}, the result follows for
any $G$-embedding $X$ and any generalised projected Richardson
variety.
\end{proof}

\subsection{Normality of projected generalised Richardson varieties}

\begin{lemma}
Let $X$ be a toroidal $G$-embedding and let $X_{u,v}^{w,x}(\Omega)$ be a
generalised Richardson variety. Let $\cL$ be a globally generated line
bundle on $X$. Then the map $H^0(X,\cL)\to H^0(X_{u,v}^{w,x}(\Omega),\cL)$
is surjective and the groups $H^i(X_{u,v}^{w,x}(\Omega),\cL)$ vanish
for $i>0$. 
\end{lemma}

\begin{proof}
We may assume that $X$ is projective. Let $X_{u,v}(\Omega)$ be the Schubert
variety and let $D$ be an ample $B\times B$-divisor. Then $D$ is a
union of irreducible components of $\partial X_{u,v}(\Omega)$ the union of
proper generalised Schubert subvarieties in $X_{u,v}(\Omega)$, it does not
contain $X_{u,v}^{w,x}(\Omega)$ and is compatibly split. In particular
$X_{u,v}(\Omega)$ is $(p-1)D$-split compatibly splitting 
$X_{u,v}^{w,x}(\Omega)$. By \cite[Theorem 1.4.8]{BK}, we get that
the map in cohomology  $H^0(X_{u,v}(\Omega),\cL)\to
H^0(X_{u,v}^{w,x}(\Omega),\cL)$ is surjective and the cohomology
groups $H^i(X_{u,v}^{w,x}(\Omega),\cL)$ vanish for $i>0$. By
\cite[Corollary 8.5]{he-thomsen}, we have that the map $H^0(X,\cL)\to
H^0(X_{u,v}(\Omega),\cL)$ is surjective concluding the proof.
\end{proof}

\begin{cor}
\label{cor-norm-gen}
The projected generalised Richardson varieties are normal.
\end{cor}

\begin{proof}
Let $\varphi:\Xt\to X$ be a morphism of $G$-embeddings with $\Xt$
smooth and toroidal. It suffices to prove that the map
$\varphi:\Xt_{u,v}^{w,x}(\Omega)\to \varphi(\Xt_{u,v}^{w,x}(\Omega))$ is
cohomologically trivial. Let $\cL$ be an ample line bundle on $X$. We
have the following commutative diagram 
$$\xymatrix{H^i(X,\cL)\ar[r]\ar[d] & H^i(\Xt,\varphi^*\cL) \ar[d]\\
H^i(\varphi(X_{u,v}^{w,x}(\Omega)),\cL)\ar[r] &
H^i(X_{u,v}^{w,x}(\Omega),\varphi^*\cL). \\}$$ 
The top horizontal map is an isomorphism because $X$ has rational
singularities while the right vertical map is surjective by the
previous lemma. This implies that the bottom horizontal map is
surjective (between trivial groups for $i>0$). By \cite[Lemma 3.3.3]{BK}
we get the result.
\end{proof}

\subsection{Compatibly split subvarieties}

Let $X$ be a complete $G$-embedding.

\begin{thm}
The compatible split subvarieties for the splitting obtained in
Proposition \ref{thm-HT} are the projected generalised Richardson
varieties.
\end{thm}

\begin{proof}
We use the following result of Knutson, Lam and Speyer (see
\cite[Theorem 5.3]{KLS}): Let $X$ be complete, normal and Frobenius
split and $\mathcal{Y}$ a finite collection of compatibly split
subvarieties of $X$ defining a stratification and satisfying:
\begin{itemize}
\item[1.] each closed stratum $Y\in\mathcal{Y}$ is normal
\item[2.] each open stratum
  $Y\setminus\cup_{Z\in\mathcal{Y},Z\subsetneq Y}Z$ is regular, and 
\item[3.] $\partial X=\cup_{Y\in\mathcal{Y},\codim_X Y=1}Y$ is an
  anticanonical divisor.
\end{itemize}
Then $\mathcal{Y}$ contains all the compatibly split subvarieties in
$X$ and for each $Y\in\mathcal{Y}$, the union
$\cup_{Z\in\mathcal{Y},Z\subsetneq Y}Z$ is an anticanonical divisor.

Let $\mathcal{Y}$ be the family
$(\Pi_{u,v}^{w,x}(\Omega))_{(\Omega,u,v,w,x) \in \mathfrak{R}_X}$ of
projected generalised Richardson varieties. By Proposition
\ref{prop-strat-gen} the family $\mathcal{Y}$ is a stratification. By
Corollary \ref{cor-norm-gen} projected generalised Richardson
varieties are normal and by Proposition \ref{prop-smooth-gen} the open
strata are smooth. Furthermore, the divisorial strata are the
divisorial generalised Richardson varieties \emph{i.e.} the divisors
stable under $G\times G$, $B\times B$ or $B^-\times B^-$. This is
exactly $\partial X$ and it is an anticanonical divisor by
\cite[Proposition 6.2.6]{BK}. The result follows. 
\end{proof}

\begin{remark}
Note that as a corollary of the above proof we have that any projected
generalised Richardson variety is of the form
$\Pi_{u,v}^{w,x}(\Omega)$ for $(\Omega,u,v,w,x) \in \mathfrak{R}_X$. 
\end{remark}

\begin{remark}
A non irreducible generalised Richardson varieties is not a projected
generalised Richardson variety. However its irreducible components are
projected generalised Richardson varieties. 
\end{remark}

\begin{cor}
The divisor $\displaystyle{\sum[\Pi_{u',v'}^{w',x'}(\Omega')]}$
where the sum runs over all codimension one projected generalised
Richardson subvarieties of $\Pi_{u,v}^{w,x}(\Omega)$ is an
anticanonical divisor in $\Pi_{u,v}^{w,x}(\Omega)$.
\end{cor}

\begin{proof}
Follows from the above result and \cite[Theorem 5.3]{KLS}.
\end{proof}

\providecommand{\bysame}{\leavevmode\hbox to3em{\hrulefill}\thinspace}
\providecommand{\MR}{\relax\ifhmode\unskip\space\fi MR }
\providecommand{\MRhref}[2]{%
  \href{http://www.ams.org/mathscinet-getitem?mr=#1}{#2}
}
\providecommand{\href}[2]{#2}

\end{document}